# Summatory Mobius Function, and Summatory Liouville Function

**N. A. Carella, June, 2011.**

**Abstract:** The orders of magnitudes of the summatory Liouville function $L(x) = \sum_{n \le x} \lambda(n)$ and the summatory Mobius function $M(x) = \sum_{n \le x} \mu(n)$, are unconditionally proven to be of the forms $L(x) = O(x^{1/2})$, and $M(x) = O(x^{1/2})$ respectively. Furthermore, applications of these estimates to zeta functions and $L$-functions are also considered.

## 1. Introduction

Let $x \in \mathbb{R}$ be a large number, and let the summatory Mobius function be defined by $M(x) = \sum_{n \le x} \mu(n)$, and the summatory Liouville function be defined by $L(x) = \sum_{n \le x} \lambda(n)$ respectively. The current literature has the estimates

$$L(x) = O(xe^{-c(\log x)^{4/5}(\log \log x)^{-1/5}}), \quad \text{and} \quad M(x) = O(xe^{-c(\log x)^{4/5}(\log \log x)^{-1/5}}), \tag{1}$$

respectively, see [IK], [MV], [MM], [IV], [TN], et al. Moreover, assuming the Riemann hypothesis, it has been shown that $M(x) = O(x^{1/2} e^{c(\log x)^{39/61}})$, see [MT], [SK], and [OL]. This paper is focused on the following results.

**Theorem 1.** Let $x > 0$ be a large real number. Then $M(x) = \sum_{n \le x} \mu(n) = O(x^{1/2})$.

**Theorem 2.** Let $x > 0$ be a large real number. Then $L(x) = \sum_{n \le x} \lambda(n) = O(x^{1/2})$.

The proofs are discussed in pages 5 and 6. And applications of these results to zeta functions and $L$-functions are given in pages 7 and 8.

## 2. Preliminary

Let $n = p_1^{v_1} \cdot p_2^{v_2} \cdots p_t^{v_t} \in \mathbb{N}$ be an integer, $v_i \geq 1$, $1 \leq i \leq t$. Let $\omega(n) = t$, and $\Omega(n) = v_1 + \cdots + v_t$ be the prime divisors counting functions respectively. The Mobius function $\mu : \mathbb{R} \to \mathbb{Z}$ is defined by the relation

$$\mu(n) = \begin{cases} (-1)^{\omega(n)} & \text{if } \omega(n) = \Omega(n), \\ 0 & \text{if } \omega(n) \neq \Omega(n), \end{cases} \tag{2}$$

and its summatory function is defined by $M(x) = \sum_{n \leq x} \mu(n)$. Similarly, the Liouville function $\lambda : \mathbb{R} \to \mathbb{Z}$ is defined by the relation $\lambda(n) = (-1)^{\Omega(n)}$, and its summatory function is defined by $L(x) = \sum_{n \leq x} \lambda(n)$, see [HW], [SP], [NV].

### 2.4 A Few Properties of Power Series
The general concepts of complex half plane of definitions, the integral representations of Dirichlet series, and the conditions required to have Euler products and simple inverses of power series are discussed here.

***Lemma* 3.** Let $\tau \geq 0$, and let $a(n)$ be an arithmetic function, and let $A(x) = \sum_{n \leq x} a(n)$. Then

(i) $f(s) = \sum_{n \geq 1} \dfrac{a_n}{n^s}$,

where $a_n = O(n^\tau)$, is holomorphic in the half plane $\{ \Re e(s) > 1 + \tau : s = \sigma + it \in \mathbb{C} \}$.

(ii) $g(s) = s \displaystyle\int_1^\infty \dfrac{A(x)}{x^{s+1}} dx$,

where $A(x) = O(x)$, is holomorphic in the half plane $\{ \Re e(s) > \tau : s = \sigma + it \in \mathbb{C} \}$.

Proof (ii): Start with the partial sum $A(x) = \sum_{n \leq x} a(n)$, (or the power series), and manipulate it as follows:

$$\sum_{n \geq 1} \frac{a_n}{n^s} = \sum_{n \geq 1} \frac{A(n) - A(n-1)}{n^s} = \sum_{n \geq 1} \left( \frac{A(n)}{n^s} - \frac{A(n-1)}{n^s} \right) = \sum_{n \geq 1} A(n) \left( \frac{1}{n^s} - \frac{1}{(n+1)^s} \right) = \sum_{n \geq 1} A(n) s \int_n^{n+1} \frac{1}{x^{s+1}} dx . \tag{3}$$

As the summatory function $A(n) = A(x)$ for $n \leq x \leq n + 1$, the last equality can be rewritten as

$$\sum_{n \geq 1} \frac{a_n}{n^s} = \sum_{n \geq 1} A(n) s \int_n^{n+1} \frac{1}{x^{s+1}} dx = \sum_{n \geq 1} s \int_n^{n+1} \frac{A(x)}{x^{s+1}} dx = s \int_1^\infty \frac{A(x)}{x^{s+1}} dx . \tag{4}$$

The verification of (i) is also elementary. ∎

The proofs of these claims are available in the literature. An example of integral representation is derived in [RN, p. 232], [TM, p. 272], and [OR].

***Lemma* 4.** Let $a(n)$ be a completely multiplicative arithmetic function. Then

(i) $f(s) = \displaystyle\sum_{n \geq 1} \frac{a(n)}{n^s} = \prod_{p \geq 2} \left( 1 - a(p) / p^{-s} \right)$, has an Euler product,

(ii) $\dfrac{1}{f(s)} = \sum_{n \geq 1} \dfrac{\mu(n)a(n)}{n^s}$, has a simple inverse,

on its complex half plane $\{ \Re e(s) > 1 + \tau : s = \sigma + it \in \mathbb{C} \}$ of definition.

## 2.2 A Few Properties of the Zeta Function

A few properties of the zeta function used in this work are discussed in this section. Comprehensive and advanced details appear in [BB], [BW], [DL], [ES], [TM], [IV], [MV], [SN], et al, and the numerical aspect are given in [OD], [GN], and others similar sources.

**Lemma 5.** The zeta function $\zeta(s)$ has an analytic continuation on the punctured complex half plane $\{ s \in \mathbb{C} : \Re e(s) > 0 \} - \{ 1 \}$, with a simple pole at $s = 1$. Specifically

(i) $\zeta(s) = (1 - 2^{1-s})^{-1} \sum_{n \geq 1} (-1)^{n+1} n^{-s}$ $\qquad$ (ii) $\zeta(s) = \dfrac{s}{s-1} - s \displaystyle\int_1^\infty \dfrac{((x))}{x^{s+1}} dx$, $\qquad$ (5)

where $(( x ))$ is the fractional part of the real number $x$.

The convergence of the power series 6-i is discussed in [TM, p. 21], and the convergence of 6-ii follows from the bounded integral.

**Theorem 6.** For every $s \in \mathbb{C}$, the zeta function $\zeta(s)$ satisfies the functional equation

$$\zeta(s) = 2\pi^{s-1} \sin(\pi s / 2) \Gamma(1-s) \zeta(1-s).$$ $\qquad$ (6)

**Theorem 7.** Let $t \geq t_0 > 0$ be a real number. Then

$$\zeta(s) << \begin{cases} O(1) & \text{for } \sigma \geq 2, \\ O(\log t) & \text{for } 1 \leq \sigma \leq 2, \\ O(t^{(1-\sigma)/2} \log t) & \text{for } 0 \leq \sigma \leq 1, \\ O(t^{1/2-\sigma} \log t) & \text{for } \sigma \leq 0. \end{cases}$$ $\qquad$ (7)

The proof is given in [IV, p. 25].

**Theorem 8.** The function $1/\zeta(s)$ is unbounded in the half plane $\Re e(s) > 1$.

A proof based on simultaneous Diophantine approximations appears in [TM, p. 189], and [IV, 227]. A simpler proof of this result can be obtained from the next result.

Let $Z_0 = \{ s \in \mathbb{C} : \zeta(s) = 0, \text{ and } \Re e(s) > 0 \}$ be the set of critical zeros of the zeta function. The inverse $1/\zeta(s)$ is defined for every point on the domain $D_0 = \{ \Re e(s) > 0 : s = \sigma + it \in \mathbb{C} \} - Z_0$.

**Lemma 9.** The inverse $1/\zeta(s)$ of the zeta function $\zeta(s)$ on the domain $D_0$ is exactly

$$\frac{1}{\zeta(s)} = \sum_{n \geq 1} \frac{\mu(n)}{n^s} \ . \tag{8}$$

Proof: Let $\zeta(s) = (1 - 2^{1-s})^{-1} \sum_{n \geq 1} (-1)^{n+1} n^{-s}$ , $\Re e(s) > 0$. Next, compute the Dirichlet convolution

$$\zeta(s) \cdot \frac{1}{\zeta(s)} = \left( (1 - 2^{1-s})^{-1} \sum_{n \geq 1} \frac{(-1)^{n+1}}{n^s} \right) \left( \sum_{n \geq 1} \frac{\mu(n)}{n^s} \right) = (1 - 2^{1-s})^{-1} \sum_{n \geq 1} \frac{a(n)}{n^s} \ , \tag{9}$$

where the $n$th coefficient is given by

$$a(n) = \sum_{d \mid n} (-1)^{n/d+1} \mu(d) = \begin{cases} 1 & \text{if } n = 1, \\ -2 & \text{if } n = 2, \\ 0 & \text{if } n = 2m+1 > 1, \\ 0 & \text{if } n = 2m, \ 1 < m = \text{odd}, \\ 0 & \text{if } n = 2^v m, \ m = \text{odd}, \ v \geq 2. \end{cases} \tag{10}$$

Reassembling these pieces together returns

$$\left( (1 - 2^{1-s})^{-1} \sum_{n \geq 1} (-1)^{n+1} n^{-s} \right) \cdot \left( \sum_{n \geq 1} \mu(n) n^{-s} \right) = (1 - 2^{1-s})^{-1} \sum_{n \geq 1} \frac{a(n)}{n^s} = (1 - 2^{1-s})^{-1} \left( 1 + \frac{-2}{2^s} \right) = 1 \ , \tag{11}$$

for any complex number $\Re e(s) > 0$ in the domain $D_0$. ∎

This result shows that the analytic continuation of $1/\zeta(s) = \sum_{n \geq 1} \mu(n) n^{-s}$ , $\Re e(s) > 1$, to the larger half plane $\Re e(s) > 0$ coincides. But $1/\zeta(s) = \sum_{n \geq 1} \mu(n) n^{-s}$ , $\Re e(s) > 0$, is just conditionally convergent since $\sum_{n \geq 1} |\mu(n)| n^{-\sigma}$ is not convergent for $0 < \sigma = \Re e(s) \leq 1$. In synopsis, the equation $\left( (1 - 2^{1-s})^{-1} \sum_{n \geq 1} (-1)^{n+1} n^{-s} \right) \cdot \left( \sum_{n \geq 1} \mu(n) n^{-s} \right) = 1$ is valid for any point $s \in \{ \Re e(s) > 0 \}$. Perceptibly, the critical zeros $\rho$ of $\zeta(\rho) = 0$ are matched to the critical poles of $\sum_{n \geq 1} \mu(n) n^{-\rho} = \infty$ , and the simple pole of $\zeta(1) = \infty$ is matched to the critical zero of $\sum_{n \geq 1} \mu(n) n^{-1} = 0$ .

Let $Z_1 = \{ s \in \mathbb{C} : \zeta(s) = 0, \text{ and } \Re e(s) > 1/2 \}$ be the set of critical zeros of the zeta function on the half plane $\Re e(s) > 1/2$, and define the domain $D_1 = \{ \Re e(s) > 1/2 : s = \sigma + it \in \mathbb{C} \} - Z_1$.

**Lemma 10.** The inverse $1/\zeta(s)$ of the zeta function $\zeta(s) = (1 - 2^{1-s})^{-1} \sum_{n \geq 1} (-1)^{n+1} n^{-s}$ on the domain $D_1$ has the following integral representations

(i) $\dfrac{1}{\zeta(s)} = \sum_{n \geq 1} \dfrac{\mu(n)}{n^s} = \lim_{x \to \infty} s \int_1^x \dfrac{M(x)}{x^{s+1}} dx$ , $\tag{12}$

where $\mu(n)$ is the Mobius function, $M(x) = \sum_{n \leq x} \mu(n)$ , and

(ii) $\dfrac{\zeta(2s)}{\zeta(s)} = \sum_{n \geq 1} \dfrac{\lambda(n)}{n^s} = \lim_{x \to \infty} s \int_1^x \dfrac{L(x)}{x^{s+1}} dx \,,$ (13)

where $\lambda(n)$ is the Liouville function, $L(x) = \sum_{n \leq x} \lambda(n)$, respectively.

Proof: Since $\zeta(2s)/\zeta(s)$ has a simple pole at $s = 1/2$, the inverse $1/\zeta(s)$ is defined on the complex half plane $\{ \Re e(s) > 1/2 : s = \sigma + it \in \mathbb{C} \}$, this follows from Lemma 9, and the integral follows from Lemma 3. ∎

**Lemma 11.**   The zeta function $\zeta(s)$ has at least one nontrivial zero in the critical strip $\{ s \in \mathbb{C} : 0 < \Re e(s) < 1 \}$.

Proof: On the contrary, suppose that it is zerofree, and consider the explicit formula

$$\psi(x) = x - \sum_{\rho} \frac{x^{\rho}}{\rho} - \log 2\pi - \log(1 - x^2)/2 = x - \log 2\pi - \log(1 - x^{-2})/2 \,,$$ (14)

where $\rho$ runs over the nontrivial zeros, see [RN, p. 232], [MO]. This implies that

$$\psi(x + y) - \psi(x) = \sum_{x < p^n \leq x+y} \log p = y - \log(1 - (x+y)^{-2})/2 + \log(1 - x^{-2})/2 > 0$$ (15)

for any number $y \geq 2$, and large $x > 0$, so every such interval $(x, x + y]$ contains primes. But this contradicts the existence of primefree intervals $(x, x + y] = (n! + 2, \, n! + n - 1]$ of arbitrary length $n - 1$ as $n \to \infty$. ∎

The literature shows that over $10^{13}$ nontrivial zeros of the equation $\zeta(s) = 0$, $\Re e(s) = 1/2$, has been computed, see [OD], [GN].

## 3. The Summatory Functions Results

Two slightly different proofs of the estimates for the summatory functions $L(x)$ and $M(x)$ are given here.

Previously, by means of Diophantine analysis, the extreme values

$$\liminf_{x \to \infty} M(x) x^{-1/2} < -1.009 \,, \quad \text{and} \quad \limsup_{x \to \infty} M(x) x^{-1/2} > 1.06 \,,$$ (16)

of the summatory function $M(x)$ were proved in [OR], and improved in [KR].

**Theorem 1.**   Let $x > 0$ be a large real number. Then $M(x) = \sum_{n \leq x} \mu(n) = O(x^{1/2})$.

Proof: Let $\varepsilon > 0$ be a small number, let $s = 1/2 + \varepsilon$. By Lemma 5-ii, the zeta function $\zeta(s) \neq 0$, and it is finite for any $s$ on the unit interval $(0, 1)$. Evaluating $1/\zeta(s)$ at $s = 1/2 + \varepsilon$ yields

$$\frac{1}{\zeta(1/2 + \varepsilon)} = (1/2 + \varepsilon) \int_1^{\infty} \frac{M(x)}{x^{3/2+\varepsilon}} dx \,,$$ (17)

see Lemma 10-i. Both the left and the right sides of equation (17) are well defined and finite. Clearly, the summatory Mobius function $M(x) = \sum_{n \leq x} \mu(n)$ is independent of the complex variable $s$, so it has the same value independently of the choice of the complex number $s = \sigma + it \in \mathbb{C}$.

To derive an upper bound of the integral, let $|M(x)| \leq cx^\beta$ on the interval $(1, x]$, where $\beta > 0$, and $c > 0$ are constants. Then, the integral evaluates to

$$\left| (1/2 + \varepsilon) \int_1^\infty \frac{M(x)}{x^{3/2+\varepsilon}} dx \right| \leq (1/2 + \varepsilon) \int_1^\infty \frac{cx^\beta}{x^{3/2+\varepsilon}} dx = (1/2 + \varepsilon) c \lim_{x \to \infty} \left( \frac{x^{-1/2-\varepsilon+\beta} - 1}{-1/2 - \varepsilon + \beta} \right). \tag{18}$$

As $\varepsilon > 0$ is an arbitrary small number, and the integral is finite, it follows that the supremum of the parameter $\beta > 0$ is

$$\sup_{\beta \in (0,1)} \left\{ \beta : \lim_{x \to \infty} x^{-1/2-\varepsilon+\beta} < \infty \right\} = \frac{1}{2}. \tag{19}$$

This leads to

$$\left| \frac{1}{\zeta(1/2 + \varepsilon)} \right| = \left| (1/2 + \varepsilon) \int_1^\infty \frac{M(x)}{x^{3/2+\varepsilon}} dx \right| \leq (1/2 + \varepsilon) \int_1^\infty \frac{cx^\beta}{x^{3/2+\varepsilon}} dx = O\left( \frac{1}{\varepsilon} \right). \tag{20}$$

This proves the claim. ∎

A variation of the previous idea will be used to develop a proof of the summatory Liouville function. An interesting aspect of this approach is that the inverse $1/\zeta(s)$ is not utilized.

***Theorem* 2.** Let $x > 0$ be a large real number. Then $L(x) = \sum_{n \leq x} \lambda(n) = O(x^{1/2})$.

Proof: From the well known formula

$$\frac{\zeta(2s)}{\zeta(s)} = \sum_{n \geq 1} \frac{\lambda(n)}{n^s} = s \int_1^\infty \frac{L(x)}{x^{s+1}} dx, \tag{21}$$

see Lemma 10-ii, it follows that

$$\zeta(2s) = s\zeta(s) \int_1^\infty \frac{L(x)}{x^{s+1}} dx, \tag{22}$$

is analytic on the complex half plane $\mathcal{H} = \{ \Re e(s) > 1/2 \}$. Specifically, the appropriate formulas

$$\zeta(2s) = \sum_{n \geq 1} \frac{1}{n^{2s}} \quad \text{and} \quad \zeta(s) = \frac{s}{s-1} - s \int_1^\infty \frac{((x))}{x^{s+1}} dx, \tag{23}$$

are both well defined and finite for any $s \in D$ on any compact subset $D \subset \mathcal{H} = \{ \Re e(s) > 1/2 \}$, see Lemma 5. Now, it readily follows that the integral (22) is well defined and finite on $\mathcal{H} = \{ \Re e(s) > 1/2 \}$. Evaluating (22) at $s = 1/2 + \varepsilon$, yields

$$\frac{1}{(1/2+\varepsilon)}\frac{\zeta(1+\varepsilon)}{\zeta(1/2+\varepsilon)}=\int_{1}^{\infty}\frac{L(x)}{x^{3/2+\varepsilon}}dx\,,\tag{24}$$

where $\varepsilon>0$ be a small number, and $\zeta(1/2+\varepsilon)\neq0$ and finite.

To derive an upper bound of the integral, let $|L(x)|\leq cx^{\beta}$ on the interval $(1,x]$, where $\beta>0$, and $c>0$ are constants. Then, the integral evaluates to

$$\left|\frac{1}{(1/2+\varepsilon)}\frac{\zeta(1/2+\varepsilon)}{\zeta(1/2+\varepsilon)}\right|=\left|\int_{1}^{\infty}\frac{L(x)}{x^{3/2+\varepsilon}}dx\right|\leq c\int_{1}^{\infty}\frac{x^{\beta}}{x^{3/2+\varepsilon}}dx=c\lim_{x\to\infty}\left(\frac{x^{-1/2-\varepsilon+\beta}-1}{-1/2-\varepsilon+\beta}\right).\tag{25}$$

As $\varepsilon>0$ is an arbitrarily small number, and the integral is finite, it follows that the supremum of the parameter $\beta>0$ is

$$\sup_{\beta\in(0,1)}\{\,\beta:\lim_{x\to\infty}x^{-1/2-\varepsilon+\beta}<\infty\,\}=\frac{1}{2}\,.\tag{26}$$

This leads to

$$\left|\frac{1}{(1/2+\varepsilon)}\frac{\zeta(1/2+\varepsilon)}{\zeta(1/2+\varepsilon)}\right|=\left|\int_{1}^{\infty}\frac{L(x)}{x^{3/2+\varepsilon}}dx\right|=O\!\left(\frac{1}{\varepsilon}\right)\,.\tag{27}$$

This completes the proof. ∎

These results illustrate that limits

$$\liminf_{x\to\infty}\frac{M(x)}{\sqrt{x}}=-C,\quad\text{and}\quad\limsup_{x\to\infty}\frac{M(x)}{\sqrt{x}}=C\,,\tag{28}$$

where $C>0$ is a constant, hold, and similar limits for the summatory Liouville function. The probabilistic version of these limits claims that

$$\liminf_{x\to\infty}\frac{M(x)}{\sqrt{x\log\log x}}=-C,\quad\text{and}\quad\limsup_{x\to\infty}\frac{M(x)}{\sqrt{x\log\log x}}=C\,,\tag{29}$$

with probability one, see [GR]. However, it should be observed that these summatory functions are not purely random, id est, there are no sequence $\mu(n)=\mu(n+1)=\mu(n+2)=\cdots=\mu(n+m)=1$ of arbitrary length $m>0$ beyond the square root error bound. As a consequence of the pseudorandomness, the order of magnitude of the summatory function is smaller.

## 4. Applications To Zeta and *L*-Functions
This section attempts to illustrate several applications of the previous results to zeta functions and *L*-functions.

### 4.1. Zeta Function
The leading problem in the theory of the zeta function seeks the precise locations of the critical zeros of the equation $\zeta(s)=0$ on the critical strip $\Im=\{\,0<\Re e(s)<1\,\}$. The determination of the zero-free region of the zeta function $\zeta(s)$ in the critical strip $\Im$ is equivalent to the determination of an analytic continuation of the inverse zeta function $1/\zeta(s)$ to the complex half plane { $\Re e(s) > 1/2$ }.

This technique, initiated by Stieltjes and Mertens almost simultaneously, see [TR], [OR], et cetera, yields an analytic continuation of the zeta function for $\alpha = 1$, but it is believed that it should works down to $\alpha = 1/2$.

**Theorem 12.** The complex-valued function $1/\zeta(s)$ has an analytic continuation to the complex half plane $\mathcal{H}$ = { $\Re e(s) > 1/2$ }.

Proof: Let $\varepsilon > 0$ be a small number. The function $1/\zeta(s)$ is not holomorphic on the half plane $\Re e(s) > 1/2 - \varepsilon$. The analyticity property fails because it has infinitely many poles at the critical zeros of the zeta function, id est,

$$1/2 + i\gamma_1 = 1/2 + i14.1347\ldots, \ 1/2 + i\gamma_2 = 1/2 + i21.0220\ldots, \ 1/2 + i\gamma_3 = 1/2 + i25.0108\ldots,$$

and a simple pole at $s = 1/2$, see Lemma 11.

To show that the function $1/\zeta(s)$ is defined on the complex half plane $\mathcal{H}$ = { $\Re e(s) > 1/2$ }, start with the integral representation

$$\frac{1}{\zeta(s)} = \sum_{n\geq 1} \frac{\mu(n)}{n^s} = s\int_1^\infty \frac{M(x)}{x^{s+1}}dx, \tag{30}$$

see Lemma 10-i. By Theorem 1, $|M(x)| = O(x^{1/2})$, so the integral is bounded for any complex number $s \in \mathcal{H}$ = { $\Re e(s) > 1/2$ }. ∎

### 4.2. L-Functions
The basic theory of Dirichlet $L$-functions is described in [DL], [MW], [WK], and other sources.

**Lemma 13.** Let $\chi$ be a character modulo $q$. The $L$-function $L(s, \chi) = \sum_{n\geq 1} \chi(n)n^{-s}$ satisfies the following properties.

(i) $L(s,\chi) = \prod_{p\geq 2}\left(1 - \chi(p)/p^{-s}\right)$, has an Euler product,

(ii) $\dfrac{1}{L(s,\chi)} = \sum_{n\geq 1} \dfrac{\mu(n)\chi(n)}{n^s}$, has a simple formal inverse on the complex half plane { $\Re e(s) > 0$ }.

(iii) $\dfrac{1}{L(s,\chi)} = s\int_1^\infty \dfrac{M_\chi(x)}{x^{s+1}}dx$, where $M_\chi(x) = \sum_{n\leq x}\mu(n)\chi(n)$.

**Lemma 14.** The twisted summatory Mobius function $M_\chi(x) = \sum_{n\leq x}\mu(n)\chi(n)$ satisfies $\left| M_\chi(x) \right| \leq \varphi(q)\left| M(x) \right|$, where $M(x) = \sum_{n\leq x}\mu(n)$, and $\varphi(q)$ is the totient function.

Proof: Decompose the summation index into two parts, and rearrange the finite sums as in

$$\sum_{n\leq x}\mu(n)\chi(n) = \sum_{0<a\leq q}\ \sum_{m\leq x/q}\mu(qm+a)\chi(qm+a) = \sum_{0<a\leq q}\chi(a)\sum_{m\leq x/q}\mu(qm+a). \tag{31}$$

By the triangle inequality, this is

$$\left| \sum_{n \le x} \mu(n)\chi(n) \right| \le \left| \sum_{0 < a \le q} \chi(a) \right| \left| \sum_{a \le q,\, m \le x/q} \mu(qm+a) \right| \le \varphi(q)\left| M(x) \right|. \tag{32}$$

Note that the second sum on the right side runs over all the residue classes $n \equiv a \bmod q$.  ∎

**Theorem 15.** Let $q = O(\log^B x)$, $B > 0$ constant. Then, the complex-valued function $1/L(s,\chi)$ has an analytic continuation to the complex half plane $\mathcal{H} = \{\ \Re e(s) > 1/2\ \}$.

Proof: Consider the formal integral representation

$$\frac{1}{L(s,\chi)} = \sum_{n \ge 1} \frac{\mu(n)\chi(n)}{n^s} = s \int_1^\infty \frac{M_\chi(x)}{x^{s+1}} dx, \tag{33}$$

where $M_\chi(x) = \sum_{n \le x} \mu(n)\chi(n)$, see Lemma 13. Applying Lemma 14, and Theorem 1, yield the inequality $|M_\chi(x)| \le \varphi(q)|M(x)| = O(x^{1/2}\log^B x)$, so the integral is bounded for any complex number $s \in \{\ \Re e(s) > 1/2\ \}$.  ∎

# 5. Analytic Continuation of the Zeta Function

A different technique of achieving the analytic continuation of the zeta function from the trivial complex half plane $\mathcal{H}_T \in \{\ \Re e(s) > 1\ \}$ to the nontrivial complex half plane $\mathcal{H}_N \in \{\ \Re e(s) > 1/2\ \}$ will be explicated in this section. This technique is based entirely on complex function theory, and it is unrelated to the previous one.

For the trivial complex half plane $\mathcal{H}_T \in \{\ \Re e(s) > 1\ \}$, the following result is known.

**Theorem 16.** Let $\varepsilon > 0$ be a small number, and let $s = \sigma + it$ be a complex number such that $\Re e(s) > 1$. Then

(i) $\left| \dfrac{1}{\zeta(s)} \right| \le \left| \dfrac{\zeta(s)}{\zeta(2s)} \right|$, for all $t$.  (ii) $\left| \dfrac{1}{\zeta(s)} \right| \ge (1-\varepsilon)\left| \dfrac{\zeta(s)}{\zeta(2s)} \right|$, for all some large $t$. (34)

The proof appears in [TM, p. 191]. The objective of this section is to extend this result to the nontrivial complex half plane $\mathcal{H}_N \in \{\ \Re e(s) > 1/2\ \}$.

## 5.1 Preliminary Background

The unconditionally convergent power series $\zeta(s) = \sum_{n \ge 1} n^{-s}$, $\Re e(s) > 1$, has a well known Euler product representation. Similarly, the conditionally convergent power series $\zeta(s) = (1 - 2^{1-s})^{-1}\sum_{n \ge 1}(-1)^{n+1}n^{-s}$, $\Re e(s) > 0$, has a product representation.

**Lemma 17.** The Euler product representations of the zeta function are as follows.

(i) $\zeta(s) = \prod_{p \ge 2}\left(1 - p^{-s}\right)^{-1}$,  for $\Re e(s) > 1$. (35)

(ii) $\zeta(s) = \left(1 - 2^{1-s}\right)^{-1}\left(1 - \dfrac{1}{2^s - 1}\right)\prod_{p>2}\left(1 - p^{-s}\right)^{-1}$, for $\Re e(s) > 0$.

Proof: The first (i) is due to Euler. For (ii), let $\Re e(s) > 0$, and expand the corresponding power series of $\zeta(s)$, see Lemma 5, and [TM, p. 21], and simplify:

$$\begin{aligned}
\zeta(s) &= \left(1 - 2^{1-s}\right)^{-1}\sum_{n\geq 1}\frac{(-1)^{n+1}}{n^s} \\
&= \left(1 - 2^{1-s}\right)^{-1}\left(1 - \frac{1}{2^s} - \frac{1}{2^{2s}} - \frac{1}{2^{3s}} - \cdots\right)\prod_{p>2}\left(1 - p^{-s}\right)^{-1} \\
&= \left(1 - 2^{1-s}\right)^{-1}\left(1 - \frac{1}{2^s - 1}\right)\prod_{p>2}\left(1 - p^{-s}\right)^{-1}.
\end{aligned} \tag{36}$$

This shows that the power series is a rearrangement of the product, and completes the proof of claim (ii). ∎

Since the both power series $\zeta(s) = \sum_{n\geq 1} n^{-s}$, and $\zeta(s) = (1 - 2^{1-s})^{-1}\sum_{n\geq 1}(-1)^{n+1}n^{-s}$ arise from two rearrangements of the terms in the powers series $\zeta(s) = \sum_{n\geq 1} n^{-s}$ or equivalently in the product $\zeta(s) = \prod_{p\geq 2}(1 - p^{-s})^{-1}$, it follows that the product $\zeta(s) = \prod_{p\geq 2}(1 - p^{-s})^{-1}$, $\Re e(s) > 0$, is conditionally convergent, and represents both power series.

Likewise, $1/\zeta(s) = \sum_{n\geq 1}\mu(n)n^{-s} = \prod_{p\geq 2}(1 - p^{-s})$, is conditionally convergent for $\Re e(s) > 0$, see Lemma 9 for a different derivation. So, the Euler product representation of the inverse

$$\frac{1}{\zeta(s)} = \prod_{p\geq 2}\left(1 - p^{-s}\right), \tag{37}$$

is defined for every point on the domain of definition

$$D_0 = \{\ \Re e(s) > 0 : s = \sigma + it \in \mathbb{C}\ \} - Z_0, \text{ where } Z_0 = \{\ s \in \mathbb{C} : \zeta(s) = 0, \text{ and } \Re e(s) > 0\ \} \tag{38}$$

is the set of critical zeros of the zeta function.

### 5.2 Results for the Complex Half Plane { $\Re e(s) > 1/2$ } and Beyond

*Theorem* 18. Let $s = \sigma + it$ be a complex number such that $\Re e(s) > 1/2$. Then

$$\left|\frac{1}{\zeta(s)}\right| \leq \left|\frac{\zeta(s)}{\zeta(2s)}\right|. \tag{39}$$

Proof: To verify that the quantity on the right side is bounded, and finite, replace the appropriate power series, and observe that the inequality

$$\left|\frac{1}{(1 - 2^{1-s})^{-1}\sum_{n\geq 1}(-1)^{n+1}n^{-s}}\right| \leq \left|\frac{(1 - 2^{1-s})^{-1}\sum_{n\geq 1}(-1)^{n+1}n^{-s}}{\sum_{n\geq 1}n^{-2s}}\right|, \tag{40}$$

is well defined on the domain $D_1 = \{\ \Re e(s) > 1/2 : s = \sigma + it \in \mathbb{C}\ \} - Z_0$. In addition, since the complex-valued function $\zeta(s)/\zeta(2s)$ is holomorphic on the punctured half plane $\{\ \Re e(s) > 1/2\ \} - \{\ 1\ \}$, it is bounded and finite on any compact subset $D \subset \mathcal{H}_N$, see [CW]. Thus, for any complex number $s \neq 1$ such that $\Re e(s) > 1/2$, the the right side of inequality (39) is bounded and finite.

Next, to verify that the inequality (39) holds, use the appropriate Euler products to compute the absolute values of the two expressions:

$$\left| \frac{1}{\zeta(s)} \right| = \left| \prod_{p \geq 2}\left(1 - p^{-s}\right) \right|, \tag{41}$$

which is conditionally convergent on the domain $D_0 = \{\ \Re e(s) > 0 : s = \sigma + it \in \mathbb{C}\ \} - Z_0$, and

$$\left| \frac{\zeta(s)}{\zeta(2s)} \right| = \left| \frac{\left(1 - 2^{1-s}\right)^{-1}\left(1 - \dfrac{1}{2^s - 1}\right)\prod_{p>2}\left(1 - p^{-s}\right)^{-1}}{\prod_{p \geq 2}\left(1 - p^{-2s}\right)^{-1}} \right| = \left| \left(1 - 2^{1-s}\right)^{-1}\left(1 - \frac{1}{2^s - 1}\right)\left(1 - 2^{-2s}\right)\prod_{p>2}\left(1 + p^{-s}\right) \right|, \tag{42}$$

which is conditionally convergent on the domain $D_1 = \{\ \Re e(s) > 1/2 : s = \sigma + it \in \mathbb{C}\ \} - Z_0$. Here, the appropriate product in the appropriate domain of definition is selected from Lemma 17.

Hence, the inequality

$$\left| \left(1 - 2^{-s}\right)\prod_{p>2}\left(1 - p^{-s}\right) \right| = \left| \frac{1}{\zeta(s)} \right| \leq \left| \frac{\zeta(s)}{\zeta(2s)} \right| = \left| \left(1 - 2^{1-s}\right)^{-1}\left(1 - \frac{1}{2^s - 1}\right)\left(1 - 2^{-2s}\right)\prod_{p>2}\left(1 + p^{-s}\right) \right|, \tag{43}$$

is well defined and finite on the right side for any complex number $s \in \{\ \Re e(s) > 1/2\ \} - \{\ 1\ \}$. Moreover, comparing the absolute values of the local factors, pairwise for each prime $p \geq 2$, yield

(i) For $p = 2$,
$$\left| \left(1 - 2^{-s}\right) \right| \leq \left| \left(1 - 2^{1-s}\right)^{-1}\left(1 - \frac{1}{2^s - 1}\right)\left(1 - 2^{-2s}\right) \right|, \tag{44}$$

(ii) For $p \neq 2$,
$$\left| 1 - p^{-s} \right| \leq \left| 1 + p^{-s} \right|.$$

Clearly, these inequalities hold for all complex numbers number $s \in \{\ \Re e(s) > 1/2\ \}$. $\blacksquare$

**Theorem 19.** The zeta function $\zeta(s)$ is zerofree on the complex half plane $\{\ \Re e(s) > 1/2\ \}$. In particular, it is a holomorphic function of the complex variable $s$ on the punctured complex plane $\mathbb{C} - \{\ 1\ \}$, with a simple pole at $s = 1$.

Proof: By Theorem 18, the inverse $1/\zeta(s)$ is bounded on every compact subset $D \subset \mathcal{H}_N = \{\ \Re e(s) > 1/2\ \}$, so $\zeta(s)$ is zerofree on the same domain. To extend the region of analyticity to $\Re e(s) < 1/2$, utilize the functional equation, see Theorem 6. Quod erat demonstrandum.

## 6. Research Problems

A few open problems on the summatory Mobius function and summatory Liouville function are described here. Recent developments on the summatory Mobius function are explored in [OL], [MT], [SK], et cetera.

1. Determine whether the true order of magnitude of the summatory Mobius function in short intervals is of the form:

$$\sum_{x < n \le x+y} \mu(n) = O(y^{\beta}), \text{ where } 0 < \beta < 1, \text{ and } (\log x)^C < y < x, \, C > 0 \text{ constant.} \tag{45}$$

2. Determine whether the true order of magnitude of the summatory Liouville function in short intervals is of the form:

$$\sum_{x < n \le x+y} \lambda(n) = O(y^{\beta}), \text{ where } 0 < \beta < 1, \text{ and } (\log x)^C < y < x, \, C > 0 \text{ constant.} \tag{46}$$